\documentclass[twoside]{article}
\usepackage[utf8]{inputenc}
\usepackage[russian,english]{babel}
\usepackage[T1,T2A]{fontenc}
\usepackage{multicol}[2]
\usepackage[sc]{mathpazo} 
\usepackage[hmarginratio=1:1,top=32mm,columnsep=20pt]{geometry} 
\usepackage{multicol} 
\usepackage{booktabs} 
\usepackage{hyperref} 
\usepackage{paralist} 
\usepackage{amsmath}
\usepackage{amsthm}
\usepackage{abstract} 
\hypersetup{unicode,colorlinks,citecolor={black}}
\newtheorem{theorem} {\bf  { Теорема}}

\usepackage{titlesec} 
\renewcommand\thesection{\Roman{section}} 
\renewcommand\thesubsection{\Roman{subsection}} 
\titleformat{\section}[block]{\large\scshape\centering}{\thesection.}{1em}{} 
\titleformat{\subsection}[block]{\large}{\thesubsection.}{1em}{} 
\title{\vspace{-15mm}\fontsize{24pt}{10pt}\textbf{A gradient catastrophe mechanism in contexts of the phase change condition }} 
\author{
\fontsize{11pt}{4pt}\textbf{A.A.Durmagambetov}\\
\normalsize \href{mailto:aset.durmagambet@gmail.com}{aset.durmagambet@gmail.com} 
\vspace{-5mm}
}
\date{}
\begin{document}
\maketitle 
\bf{ABSTRACT }\\
The paper describes a gradient catastrophe mechanism in contexts of the phase change conditions.  It shows that classical methods of the function estimation theory are not suited to study a gradient catastrophe problem.  The paper presents data showing that embedding theorems do not allow to study a process of a catastrophe formation. In fact, the paper justifies Terence Tao’s pessimism about a failure of modern mathematics to solveg the Navier-Stokes problem.  An alternative method is proposed for studying the gradient catastrophe by studying Fourier transformation for a function and selecting a function singularity through phase singularities of Fourier transformation for a given function.
\section {INTRODUCTION}
The research presents a process of gradient catastrophe formation under conditions of phase change. The paper shows that classical methods of the function estimation theory in contexts of Sobolev- Schwartz Space Theory are not suitable for studying gradient catastrophe problem. Data presented shows that the embedding theorems do not allow to study a process of a catastrophe formation. Actually, the paper justifies Terence Tao’s pessimism about a failure of using present mathematical methods for solving the Navier-Stokes problem.  An alternative method is proposed for studying gradient catastrophe by studying Fourier transformation for a function and selecting function singularity through phase singularities of Fourier transformation for a given function. It will be recalled a general definition of a gradient catastrophe - an unbounded increase of a function derivative upon conditions of boundedness of the function itself. This phenomenon occurs in various problems of hydrodynamics, such as a formation of shock waves, weather fronts, hydraulic and seismic fracturing, and others. In modern physics and mathematics, as well as in many other areas of science and technology, this phenomenon is considered as a very difficult problem, both from a theoretical and applied perspective.  From a theoretical point of view this is important as we have to know how to describe qualitative changes in processes, which are manifested in appearance of new quality objects during a process of description model evolution, and in the context of applied research, the problem is facing numerical instability in the event of a gradient catastrophe formation. Thus, we approach an important obstacle while using modeling - a barrier created by the gradient catastrophe. Since, on the one hand, the gradient catastrophe is still unknown phenomenon,  it is very important from a practical point of view, because the phenomenon is connected with the most interesting and important aspects of reality. Terence Tao formulated and illustrated this in [1] based on the Millennium problem stated by Clay Institute for the Navier-Stokes equations. Our point of view on these issues coincides with stated in article [1], but in our research we propose approaches to solving these problems.
Our point of view is that the modern mathematical methods of the theory of functions dedicated to the function estimation have ignored such an important component of the Fourier transformation as its phase.  Our research is outlined as follows:  first, we give examples of the gradient catastrophe caused by the phase change, and then proceed to an expansion of classes of functions subjected to the gradient catastrophe. Our final results lie in the nonlinear representation of functions showing some new classification of functions through a phase classification. In addition, the notions of discreteness and continuity of functions are naturally merged. And, as we think,  this leads to understanding of how discrete objects are born under a continuous change of the world. Discrete objects are associated with discrete spectrum of the Liouville- Schrödinger equations. And they, as is known, reflect the wave nature of things.  But here, we abstract away from the quantum formalism, because our goal  lies in a purely mathematical approach to the analysis of the arbitrary functions. For the analysis of which, we formally consider a function as a potential of the Schrödinger equation.  At the same time we come across the concepts that generated by the Liouville- Schrödinger equations. These concepts allow to classify and estimate functions by a phase generated by discrete spectrum of the Liouville equation.
 
 \section{RESULTS}
 Let us consider one-dimensional function $ {f}  $ and its Fourier transformation  $  \tilde{f} $. Using  notions of module and phase, we write Fourier transformation in the following form $  \tilde{f}=|\tilde{f}|\exp(i\phi) $ , where $ \phi $ is phase.  To cite  Plancherel equality: $ ||f||_{L_{2}}=Const||\tilde{f} ||_ {L_{2}} $. Here we can see that a phase is not contributed to determination of X norm. To estimate a maximum we have a simple estimate as  $  max|f|^2 \le 2||f||_{L_{2}} ||\nabla f||_{L_{2}}   $.Now we have an estimate of the function maximum in which a phase is not involved. Let us consider a behavior of a progressing wave running with a constant velocity of $  v=a$ described by function $ {F(x,t)=f(x+at)}  $.  For its Fourier transformation along x variable we have  $  \tilde{F}=\tilde{f}exp(iatk) $.  Again in this case we can see that when we will be studying a module of the Fourier transformation, we will not obtain major physical information about the wave, such as its velocity and location of the wave crest because of $ |\tilde{F}|=|\tilde{f}| $ . These two examples show 	w  eaknesses of studying Fourier transformation.  On the other hand, many researchers focus on the study of functions using embedding theorem, but in the embedding theorems main object of the study is module of function. But as we have seen in given examples, a phase is a main physical characteristic of a process, and as we can see in the mathematical studies, which use embedding theorems with energy estimates, the phase disappears. Along with phase, all reasonable information about physical process disappears, as demonstrated by Terence Tao [1] and other research considerations.  In fact, he built progressing waves that are not followed energy estimates. Let us proceed with more essential analysis of influence of the phase on behavior of functions.
 \begin  {theorem}
 . There are functions of $ W_{2}^{1}(R)$ with a constant rate of the norm for a gradient catastrophe of which a phase change of its Fourier transformation is sufficient. Proof: To prove this, we consider a sequence of testing functions $  \tilde{f_{n}}=\Delta/(1+k^2),\Delta=(i-k)^n/(i+k)^n $. it is obvious  that  $|\tilde{f_{n}}|=1/(1+k^2).$
 $  max|f_{n}|^2 \le 2||f_{n}||_{L_{2}} ||\nabla f_{n}||_{L_{2}} \le Const   $.. Calculating the Fourier transformation of these testing functions, we obtain: $f_{n}= x(-1)^{(n-1)}2\pi \exp(-x)L^1_{(n-1)}(2x)$ where  $ L^1_{(n-1)}(2x)$ is a Laguerre polynomial. Now we see that the functions are equibounded and derivatives of these functions will grow with the growth of $ {n}.  $
 Thus, we have built an example of a sequence of the bounded functions of  $ W_{2}^{1}(R)$ which have a constant norm  $ W_{2}^{1}(R)$  and this sequence converges to a discontinuous function.
 Thus, we have demonstrated an importance of the phase and that the phase is not involved into energy norms that are inherent to the mathematical arguments used in physical processes analysis. Our next goal is to maximally expand this class of functions in which a phase is important. Our goal is also to use a phase, which appears in the inverse scattering problem; moreover we will be interested mainly in a phase generated by a discrete spectrum of the Liouville equation. Thereby, we come now to an important subject of our research, such as an occurrence of discontinuities, fronts and other instable states in numerical modeling and which are at the same very stable physical objects. As we think, our
  \end  {theorem}
  \begin  {theorem}
Theorem 1. There are functions of $ W_{2}^{1}(R)$ with a constant rate of the norm for a gradient catastrophe of which a phase change of its Fourier transformation is sufficient. Proof: To prove this, we consider a sequence of testing functions $  \tilde{f_{n}}=\Delta/(1+k^2),\Delta=(i-k)^n/(i+k)^n $. it is obvious  that  $|\tilde{f_{n}}|=1/(1+k^2).$
$  max|f_{n}|^2 \le 2||f_{n}||_{L_{2}} ||\nabla f_{n}||_{L_{2}} \le Const   $.. Calculating the Fourier transformation of these testing functions, we obtain: $f_{n}= x(-1)^{(n-1)}2\pi \exp(-x)L^1_{(n-1)}(2x)$ where  $ L^1_{(n-1)}(2x)$ is a Laguerre polynomial. Now we see that the functions are equibounded and derivatives of these functions will grow with the growth of $ {n}.  $
Thus, we have built an example of a sequence of the bounded functions of  $ W_{2}^{1}(R)$ which have a constant norm  $ W_{2}^{1}(R)$  and this sequence converges to a discontinuous function.
 \end  {theorem}
 Thus, we have demonstrated an importance of the phase and that the phase is not involved into energy norms that are inherent to the mathematical arguments used in physical processes analysis. Our next goal is to maximally expand this class of functions in which a phase is important. Our goal is also to use a phase, which appears in the inverse scattering problem; moreover we will be interested mainly in a phase generated by a discrete spectrum of the Liouville equation. Thereby, we come now to an important subject of our research, such as an occurrence of discontinuities, fronts and other instable states in numerical modeling and which are at the same very stable physical objects. As we think, our arguments are very important in issues of plasma stability in nuclear fusion technology, since the gradient catastrophe formation serves as a preamble to a process of nuclear fusion stop. To build more in-depth analysis we apply results of scattering theory to our problem. For this, we consider a spectral problem for the Liouville equations with a potential q that satisfies and belongs  to M space of functions with the following norm $||q||_{M}=\int\limits_{-\infty}^{+\infty} |q(x) |(1+|x|)dx $
 As it is known from \begin{equation}
 	-\Psi^{"}  +q\Psi  =|k|^{2}\Psi, ~k\in C  
 \end{equation} 
 with the following asymptotics:
 \begin{gather}\label{eq:eq2}
 	\lim\limits_{x \rightarrow -\infty }\Psi_{1 } (k,x)=e^{ikx}+s_{12}(k)\exp(-ikx),\ \lim\limits_{x \rightarrow +\infty }\Psi_{1 } (k,x)=s_{11}(k)\exp(ikx) \\ 
 	\lim\limits_{x \rightarrow -\infty }\Psi_{2 } (k,x)=s_{22}(k)\exp(-ikx),\,\,  \lim\limits_{x \rightarrow +\infty }\Psi_{2 } (k,x)=\exp(-ikx)+s_{11}(k)\exp(ikx)  
 \end{gather}
 It is also known from the theory of equations [2], that any solution is a combination of some fundamental solutions satisfying certain boundary conditions.
 \begin{gather}\label{eq:eq4}
 	\lim\limits_{x \rightarrow \infty } f_{+}(k,x)\exp(-ikx)  =1,\,\  \lim\limits_{x\rightarrow-\infty}f_{-}(k,x)\exp(ikx)=1. 
 \end{gather}
 It is known [2], that they satisfy the following equations:
 \begin{gather}\label{eq:eq4}
 	f_{+}(k,x) =\exp(ikx)-\int\limits_{x}^{+\infty}  \frac{\sin(k(x-t))}{k}q(t) f_{+}(k,t) dt  ,  \\
 	f_{-}(k,x)=\exp(-ikx)+\int\limits_{-\infty}^{x}  \frac{\sin(k(x-t))}{k}q(t) f_{+}(k,t) dt \\
 	f_{+}(k,x) =\exp(ikx)-\int\limits_{-\infty}^{+\infty}  G_{+}(k,x,t)q(t) f_{+}(k,t) dt  ,  \\
 	f_{-}(k,x)=\exp(-ikx)+\int\limits_{-\infty}^{+\infty} G_{-}(k,x,t)q(t) f_{+}(k,t) dt \\
 	f_{+} =E_{+}-\sum_{j=1}^{\infty}  G_{+}^{j}E_{+}  ,\,\, f_{-} =E_{-}+\sum_{j=1}^{\infty}  G_{-}^{j}E_{-} ,\,\, \\ E_{+}(k,x)= \exp(ikx)\,\,\ E_{-}(k,x)= \exp(-ikx)   \\
 \end{gather}
 Let us also provide known results for the scattering coefficients and fundamental solutions outlined in [2].
 \begin{gather}\label{eq:eq4}
 	s_{11}f_{+}(k,x) =s_{12}f_{-}(k,x) +f_{-}(-k,x),  
 	\,\,\  s_{22}f_{-}(k,x)=s_{21}f_{+}(k,x) +f_{+}(-k,x).\\
 	s_{11}s_{12}^{*}+s_{12},\,\,s_{22}^{*}=0,\,\,s_{11}^2+s_{12}^2=s_{22}^2+s_{21}^2=1,  s_{i,j}(-k)=s_{i,j}^{*}(k), \,\  \\   \lim\limits_{|k| \rightarrow \infty } s_{12}=s_{21}=1+ O(1/|k|), \,\, 
 	,\,\ \lim\limits_{|k| \rightarrow \infty }  s_{11}=s_{22}= O(1/|k|),\, \\ s_{11}=\exp( \frac{1}{2\pi i} \int\limits_{-\infty}^{+\infty} \frac{ln (1-|s_{12}|) }   {k^{'}-k}dk^{'}    
 	\prod\limits_{j=1}^{n} \Big(\frac{ iE_{j}+k }   {k-iE_{j}}\Big) dk^{'} ,\, \\
 	s_{11}(k)=\lim\limits_{\epsilon \rightarrow 0 } =s_{11}(k+i\epsilon),\,\
 	s_{21}(k) =-     \frac{s_{12}(-k)s_{11}(k) }{s_{11}(-k)}   \\
 	s_{12}(k) =\frac{ \frac{1}{2ki } \int\limits_{-\infty}^{+\infty}\exp(-ikt)q(t) f_{+}(k,t) dt }  {1-   \frac{1}{2k i } \int\limits_{-\infty}^{+\infty}\exp(ikt)q(t) f_{+}(k,t)  dt}     
 	\\
 	s_{21}(k) =\frac{ \frac{1}{2ki} \int\limits_{-\infty0}^{+\infty}\exp(ikt)q(t) f_{-}(k,t) dt }  {1-   \frac{1}{2k i } \int\limits_{-\infty}^{+\infty}\exp(ikt)q(t) f_{-}(k,t)  dt}
 \end{gather}
 Further, we will use a modification of the last expression known in the theory of scattering  [2]
 \begin{theorem}
 	Theorem 2.  For  fundamental solutions  the following equalities are true.
 	\begin{gather}
 		\label{eq4} 
 		s_{11}f_{+}(k,x) =s_{12}f_{-}(k,x) +f_{-}(-k,x),  \\
 		\,\,\  s_{22}f_{-}(k,x)=s_{21}f_{+}(k,x) +f_{+}(-k,x).
 	\end{gather}
  \end{theorem}
  Now we are able to return to our question of the gradient catastrophe for more general class of functions. For this we consider Liouville equation and a sequence of inverse scattering problems with constant in module scattering coefficients   $s_{i,j}$, where discrete eigenvalues $ E_{i},  0<i<n+1 $  such that  $ \lim\limits_{n \rightarrow \infty} =E_{\infty}  $
  Theorem 3.  There are potentials from $ W_{2}^{1}(R) M   $ with the constant norm of $ W_{2}^{1}(R) M   $  for the gradient catastrophe for which existence of limit point for the discrete spectrum with given potential in the Liouville equation is sufficient.
  Proof. Following notations [2], we introduce functions   $ A_{+},B_{+},\Omega _{+} $   according to the formulas:
  \begin{gather}
  	\label{eq4} 
  	s_{21}(k) =\int\limits_{-\infty}^{+\infty}A_{+}(t)\exp(2ikt)dt,    \,\ 
  	\Omega _{+}(t)=\sum_{i=1}^n M_{j}^{1}exp(-E_{j}t)+A_{+}(t)\\
  	B_{+}(x,y)+\int\limits_{0}^{+\infty}B_{+}(x+y+t)\Omega _{+}(x+y+t) dt +\Omega _{+}(x+y)=0 
  \end{gather}
  where  $ M_{j}^{1}     $ are normalized numbers.  
  In other words, we will consider inverse problems of the potential recovery, and for the n-th potential we will consider a case 	with an accuracy up to n discrete eigenvalues. It is sufficient to consider a first approximation of these equations. In a first approximation, the n-th potential is recovered by the equation for  $B_{+}(x,y)  $ and also in a first approximation. We have the following arguments for the first approximation
  \begin{gather}
  	\frac{d}{dx}B_{+}(x,x)= -\frac{d}{dx}\Omega _{+}(2x). \,\  \\
  	\frac{d}{dx}\Omega _{+}(2x)=\sum_{i=1}^n -E_{j}M_{j}^{1}exp(-E_{j}t)+\frac{d}{dx}A_{+}(x)
  \end{gather}
  For the last term, we also consider a first approximation
  \begin{gather}
  	\label{eq5} 
  	\frac{d}{dx}A_{+}(x)=\int\limits_{-\infty}^{+\infty} \tilde{q}_{+}(2t) \exp(2ixt)\delta^2) dt,\,\
  	\delta=\prod\limits_{j=1}^{n} \Big(\frac{ iE_{j}+k }{k-iE_{j}}\Big) *\exp( \frac{1}{2\pi i} Vp\int\limits_{-\infty}^{+\infty} \frac{ln (1-|s_{12}|) )}   {k^{'}-k}dk^{'} )    
  	\end{gather}
  		To prove this, let us consider a sequence of $\frac{d}{dx}A_{+}(x)$with  n going  to infinity and under a proper selection of  scattering coefficients, we fall into conditions of the Theorem 1.
  		Let us come down from specific obvious examples to more systematic analysis of the gradient catastrophe. In given below all our arguments will be based on well-known equation:  
  		$$s_{21}(k) =- s_{12}(-k)s_{11}(k) /s_{11}(-k)   $$  
  		Let us consider  $s_{21},s_{12}   $ in the following form: 
  		\begin{gather}
  			2iks_{12}(k)= \tilde{q(2k)}+ I_{12}(k)               \\
  			2iks_{21}(k)= \tilde{q(-2k)} + I_{21}(k)  \\
  		\end{gather}
  		Let us conceive $\tilde{q(2k)}=U+iV $. Then we will have the following equation for  $U,V $
  		\begin{gather}
  			\label{eq46} 
  			U+iV=  2iks_{12}(k)- I_{12}(k)               \\
  			U-iV= 2iks_{21}(k) - I_{21}(k)          \\  
  			\frac{s_{11}(k)}{s_{11}(-k)} =exp(2i\delta),\, \delta=arg(s_{11}(k) ),\, \phi=2\delta  
  		\end{gather}
  		Now we can formulate the following theorem. 
  \begin{theorem}
 	 		\label{thm:t7}
 			\begin{gather}
 				\label{eq461} 
 				U=\frac{(1+\cos (\phi)R_{12}+\sin(\phi) R_{21} }{\sin(\phi)}\\
 				V=\frac{(-1+sin(\phi)) R_{12} +(1-\cos (\phi))R_{21} }{\sin(\phi)} \\                       
 			\end{gather}
 		\end{theorem}
 		Proof: Using equation () and representation for Fourier transformation we obtain
 		Whence, solving the equation for U and V, we obtain 
 		\begin{gather}
 			\label{eq46} 
 			U+iV+I_{12}(k)= (U-iV + I_{21}(k))(\cos (\phi)+i\sin(\phi) )             \\
 			U(1-\cos (\phi)))+V(1-\sin(\phi))= R_{12}          \\
 			-U\sin(\phi)+V(1+\cos (\phi))=R_{21}          \\
 			R_{12}=Real\Big(-I_{12} + I_{21}\cos (\phi) +i(I_{12}\sin(\phi)  )       \Big) \\
 			R_{21}=Im \Big(-I_{12}+I_{21}\cos (\phi)  +i(I_{12}\sin(\phi)   \Big)  
 		\end{gather}
 		\begin{theorem}
 		Theorem  5. The following estimates are true for Fourier transformation
 		\begin{gather}
 			|U| \le   C(|R_{12}| +|R_{21}| + |\nabla R_{12}| +|\nabla R_{21}| )                      \\
 			|V| \le    C(|R_{12}| +|R_{21}| + |\nabla R_{12}| +|\nabla R_{21}| )                            \\
 			\tilde{|q|} \le C(|R_{12}| +|R_{21}| + |\nabla R_{12}| +|\nabla R_{21}| )   
 		\end{gather}
 	\end{theorem}
 	 	The proof follows from the representation of U, V.  Here, we just point out this as a separate theorem in order to emphasize the significance of this result.  We note separately the terms with a derivative $\nabla R_{21},\nabla R_{12}$. Obviously, these terms are appeared due to points of the phase nulling.
 	Theorem  6.  For estimation of a maximum of the potential the following estimates are true.
 	\begin{gather}
 		|q| \le C\int\limits_{-\infty}^{+\infty}( C(|I_{12}| +|I_{21}|+|\nabla I_{12}| +|\nabla I_{21}|)   )dk                     
 	\end{gather}The proof follows from the estimation of U, V and use of $R_{12}, R_{21}$  which are simple arguments. Here we outline the theorem in order to emphasize importance of this result???. Analyzing the last formula, we see again an effect of the phase on the function behavior. In addition, a  finiteness of the discrete spectrum  is the main requirement of the  gradient catastrophe nonoccurrence. And from other hand, in case of unconstrained growth of points in discrete spectrum, we fall into the terms of theorems 1 and 2.   The last theorem expands a class of functions described in in Theorem 1, as we planned at the beginning. Now, studying the behavior of a gradient depending on q we come to the conclusion that its unconstrained growth will be dictated by the phase cluster point, which, in its turn, is due to discrete spectrum acquisition.  Hence we get the most important conclusion - we get information about the catastrophe with discrete jumps!
 	Theorem 7. For a potential the following representation is true
 	\begin{theorem}\label{thm:t8}
 		Для  потенциала  справедливо  представление
 		\begin{gather}
 			q= Q(q,E_{1},...E_{n});                     
 		\end{gather}
 	\end{theorem}The proof just consists in calculating $I_{12}(k),I_{21}(k)$ in a form of series of q and substitution of a result of the calculation into the formula for U, V moreover a right side of the obtained formula contains second-order terms only. This representation, in contrast to the classical inverse problems, allows using arbitrary information on the potential for closure of these equations, because a skeleton of this integral equation is represented by sets of constants in the form of eigenvalues.  One of the surprising properties of this representation and all this research is discreteness in continuity. Since a value of the phase, as we can see, changes discontinuously, while a potential-function itself may vary continuously. This implies an important conclusion about the instability in numerical methods, i.e. it is necessary to control phase jumps in numerical modeling to avoid falling into a state of instability. A conclusion of non-scalability of such models is critically important since eigenvalues may appear or may disappear under changes in the potential scale, whereupon a model will be changed significantly. This theorem shows that we have obtained fundamentally new nonlinear integral relations that allow taking a fundamentally fresh look at the problem of estimating functions. Now, instead of integral representations, that generate embedding theorem in the Sobolev spaces and  by which numerous outstanding achievements in modern mathematics have been gained, we turn to the newest non-linear integral relations and hope thereby opening up new pages of mathematics that will take us further into the wonderful world of mathematics. We consider this work as a starting point and assume to do series of works in this direction for a complete presentation of the conclusions of this paper.
 \section{CONCLUSION}	In this paper, we have gained excellent properties for scattering coefficients and showed, based on a number of examples, what is a cause of a gradient catastrophe. Critically important is that these results allow to illustrate results of [3], [4] using elementary tools.  As we think, after this research an attitude to a phase will be changed and that gives an opportunity to more freely handle functions undergoing a gradient catastrophe, because now we have a whole arsenal of tools for studying the gradient catastrophe in terms of the Fourier transformation in which function singularities just transform into the properties of confined and differentiable function- a phase of the scattering.
 	ACKNOWLEDGEMENTS
 	We thank the Ministry of Education and Science for this grant allocated for the project of Comprehensive Analysis in the Robotic Systems, Development of Scientific and Technological Solution, and Conducting Theoretical and Experimental Studies of Robotic System Prototypes.


\begin{thebibliography}{9}
\bibitem{b16}	Terence Tao, “Finite time blowup for an averaged three-dimensional Navier-Stokes equation,” -arXiv:1402.0290 [math.AP]
\bibitem{A002}  Фаддев Л.Д.-ДАН СССР,1958, т 121,с.63.  Гл.7 \S   5; Предисловие 
\bibitem{a02}	Asset Durmagambetov, Leyla Fazilova. Global Estimation of the Cauchy Problem Solutions’ Fourier Transform Derivatives for the Navier-Stokes Equation
International Journal of Modern Nonlinear Theory and Application Vol.2 No.4, December 2, 2013
\bibitem{b15}Navier-Stokes Equations—Millennium Prize Problems
Natural Science Vol.7 No.2,  Pub. Date: February 27, 2015
\end{thebibliography}
\end{document}